\documentclass[12pt,twoside]{amsart}
\usepackage{amsfonts,amsmath}
\usepackage[all,knot,arc]{xy}

\newcommand{\ad}{{\rm ad\,}}

\newtheorem{defn}{Definition}[section]
\newtheorem{lemma}[defn]{Lemma}

\newtheorem{theorem}[defn]{Theorem}
\newtheorem{definition}[defn]{Definition}
\newtheorem{remark}[defn]{Remark}
\newtheorem{proposition}[defn]{Proposition}
\newtheorem{corollary}[defn]{Corollary}
\newtheorem{conjecture}[defn]{Conjecture}

\newtheorem{example}[defn]{Example}

\newcommand{\zz}{{\mathbb Z}}
\newcommand{\rr}{{\mathbb R}}

\newcommand{\froyshov}{Fr{\o}yshov}
\newcommand{\ozsvath}{Ozsv\'{a}th}
\newcommand{\szabo}{Szab\'{o}}
\newcommand{\spin}{\ifmmode{\rm Spin}\else{${\rm spin}$\ }\fi}
\newcommand{\spinc}{\ifmmode{{\rm Spin}^c}\else{${\rm spin}^c$\ }\fi}
\newcommand{\spinct}{\mathfrak t}
\newcommand{\spincs}{\mathfrak s}

\newcommand{\detQ}{\delta}
\newcommand{\Diag}{\Delta}
\newcommand{\rigid}{rigid\ }
\newcommand{\rigidpunc}{rigid}

\newenvironment{narrow}[2]{%
 \begin{list}{}{%
  \setlength{\topsep}{0pt}%
  \setlength{\leftmargin}{#1}%
  \setlength{\rightmargin}{#2}%
  \setlength{\listparindent}{\parindent}%
  \setlength{\itemindent}{\parindent}%
  \setlength{\parsep}{\parskip}%
 }%
\item[]}{\end{list}}

\newif\ifpic
\picfalse    
\pictrue   

\DeclareMathOperator\PSL{PSL}
\DeclareMathOperator\im{Im}
\DeclareMathOperator\re{Re}

\begin{document}

\title{A characterisation of the $\mathbf{n\langle1\rangle%
\oplus\langle3\rangle}$ form
and applications to rational homology spheres}
\author{Brendan Owens and Sa\v{s}o Strle}
\date{\today}
\thanks{S.Strle was  supported in part by the MSZS of the
Republic of Slovenia research program No. P1-0292-0101-04 and research
project No. J1-6128-0101-04.}

\begin{abstract}
We conjecture two generalisations of Elkies' theorem on unimodular quadratic forms to 
non-unimodular forms.  We give some evidence for these conjectures including a result for
determinant 3.  These conjectures, when combined with results of \froyshov~ and 
of \ozsvath~ and \szabo,
would give a simple test of whether a rational homology 3-sphere may bound a negative-definite 
four-manifold.  We verify some predictions using Donaldson's theorem.  Based on this we compute
the four-ball genus of some Montesinos knots.
\end{abstract}

\maketitle

\pagestyle{myheadings}
\markboth{BRENDAN OWENS AND SA\v{S}O STRLE}{A CHARACTERISATION OF THE
 $\,n\langle1\rangle\oplus\langle3\rangle$ FORM}

\section{Introduction}
\label{sec:intro}

Let $Y$ be a rational homology three-sphere and $X$ a smooth
negative-definite four-manifold bounded by  $Y$.  For any $\spinc$
structure $\spinct$ on $Y$ let $d(Y,\spinct)$ denote the
correction term invariant of \ozsvath~ and \szabo~ \cite{os4}. It is
shown in \cite[Theorem 9.6]{os4} that for each $\spinc$ structure
$\spincs\in\spinc(X)$,
\begin{equation}
\label{eqn:thm9.6} c_1(\spincs)^2+{\rm
rk}(H^2(X;\zz))\le4d(Y,\spincs|_Y).
\end{equation}
This is analogous to gauge-theoretic results of \froyshov.  
These theorems constrain the
possible intersection forms that $Y$ may bound.  The above
inequality is used in \cite{bs} to constrain intersection forms of
a given rank bounded by Seifert fibred spaces, with application to
four-ball genus of Montesinos links.  In this paper we attempt to
get constraints by finding a lower bound on the left-hand side of
(\ref{eqn:thm9.6}) which applies to forms of any rank. This has been
done for unimodular forms by Elkies:

\begin{theorem}[\cite{elkies}]
\label{thm:elkies}
 Let $Q$ be a negative-definite unimodular integral quadratic
form of rank $n$.  Then there exists a
characteristic vector $x$ with $Q(x,x)+n\ge0$; moreover 
the inequality is strict unless 
$Q=n\langle-1\rangle$.
\end{theorem}

Together with (\ref{eqn:thm9.6}) this implies that an integer homology
sphere $Y$ with $d(Y)<0$ cannot bound a negative-definite
four-manifold, and if $d(Y)=0$ then the only definite pairing that $Y$
may bound is the diagonal form.  Since $d(S^3)=0$ this generalises
Donaldson's theorem on intersection forms of closed four-manifolds
\cite{d}.

In Section \ref{sec:conjectures} we conjecture two generalisations of Elkies' theorem
to forms of arbitrary determinant.  We prove some special cases,
including Theorem \ref{thm:det3} which is 
a version of Theorem \ref{thm:elkies} for forms of
determinant 3. This implies the following

\begin{theorem}
\label{cor:det3}
Let $Y$ be a rational homology sphere with $H_1(Y;\zz)=\zz/3$ and let $\spinct_0$ be the 
\spin structure on $Y$.
If Y bounds a negative-definite four-manifold $X$ then either
$$d(Y,\spinct_0)\ge -\frac12,$$
or
$$\displaystyle\max_{\spinct\in\spinc(Y)}d(Y,\spinct)\ge\frac16.$$
If equality holds in both then the intersection form of $X$ is diagonal. 
\end{theorem}

In Section \ref{sec:examples} we discuss further topological implications of our
conjectures; in particular some predictions for Seifert fibred
spaces may be verified using Donaldson's theorem.
We find two families of Seifert fibred rational homology
spheres, no multiple of which can bound negative-definite manifolds.
We use these results to determine the four-ball genus
of two families of Montesinos knots, including one whose members
are algebraically slice but not slice.
\section{Conjectured generalisations of Elkies' theorem}
\label{sec:conjectures}

We begin with some notation.  A quadratic form $Q$ of rank $n$ over
the integers gives rise to a symmetric matrix with entries
$Q(e_i,e_j)$, where $\{e_i\}$ are the standard basis for $\zz^n$;
we also denote the matrix by $Q$.  
Let $Q'$ denote the induced form on the dual $\zz^n$; 
this is represented by the inverse matrix.
Two matrices $Q_1$ and $Q_2$
represent the same form if and only if $Q_1=P^T Q_2 P$ for some
$P\in GL(n,\zz)$.

We call $y\in\zz^n$ a {\it characteristic covector} for
$Q$ if
$$(y, \xi) \equiv Q(\xi,\xi) \pmod2 \quad\forall \xi\in\zz^n.$$

We call $x\in\zz^n$ a {\it characteristic vector} for $Q$
if
$$Q(x,\xi) \equiv  Q(\xi,\xi) \pmod2 \quad\forall \xi\in\zz^n.$$

Note that the form $Q$ induces an injection 
$x\mapsto Qx$ from $\zz^n$ to its dual
with the quotient group having
order $|\det Q|$; with respect to the standard bases this map is multiplication
by the matrix $Q$.  For unimodular forms this gives
a bijection between characteristic vectors and characteristic covectors;
in general not every characteristic covector is a characteristic vector.  
Also for odd
determinant, any two characteristic vectors are congruent modulo
$2$; this is no longer true for even determinant.

Let $Q$ be a negative-definite integral form
of rank $n$ and let $\detQ$ be the absolute value of its determinant.
Denote by $\Diag=\Diag_{\detQ}$ the
diagonal form $(n-1)\langle-1\rangle\oplus\langle-\detQ\rangle$.
Both of the following give restatements of Theorem
\ref{thm:elkies} when restricted to unimodular forms.

\begin{conjecture}
\label{conj:vec} Every
characteristic vector $x_0$ is congruent modulo $2$ to a
vector $x$ with 
$$Q(x,x) + n \ge 1 - \detQ;$$
moreover the inequality is strict
unless $Q=\Diag_\detQ$.
\end{conjecture}

\begin{conjecture}
\label{conj:covec} There
exists a characteristic covector $y$ with 
$$Q'(y,y) + n \ge 
\left\{\begin{array}{ll}
1 - 1/\detQ &  \mbox{if $\detQ$ is odd,}\\
1 &  \mbox{if $\detQ$  is even;}
\end{array}\right. 
$$ 
moreover the inequality is strict
unless $Q=\Diag_\detQ$.
\end{conjecture}

We will discuss the implications of these conjectures in Section \ref{sec:examples}.

\begin{proposition}
Conjecture  \ref{conj:vec} is true when restricted to forms of rank $\le3$, and Conjecture 
\ref{conj:covec} is true when restricted to forms of rank $2$ and odd determinant.
\end{proposition}
\proof
We will first establish Conjecture \ref{conj:vec} for rank 2 forms.
In fact we prove the following stronger statement:
if $Q$ is a negative-definite form of rank 2 and determinant $\detQ$, 
then for any $x_0\in\zz^2$,
\begin{equation}
\label{eqn:firstrk2}
\displaystyle\max_{x\equiv x_0 (2)}Q(x,x)\ge -1-\detQ,
\end{equation}
and the inequality is strict unless $Q=\Diag$.

Every negative-definite
rank 2 form is represented by a {\it reduced} matrix
$$Q=\left(\begin{matrix} a & b\\ b & c \end{matrix}\right),$$
with $0\ge2b\ge a\ge c$ and $-1\ge a$.
Any vector $x_0$ is congruent modulo 2
to one of  $(0,0), (1,0), (0,1), (1,-1)$; all of these satisy 
$x^T Q x \ge a+c-2b$. Thus it suffices to show 
\begin{equation}
\label{eqn:secondrk2}
a+c-2b\ge -1-\detQ.  
\end{equation}
Note that equality holds in (\ref{eqn:secondrk2}) if $Q=\Diag$.  Suppose now that $Q\ne \Diag$.
Let
$Q_\tau=\left(\begin{matrix} a+2\tau & b+\tau\\ b+\tau & c \end{matrix}\right),$
and let $\detQ_\tau=\det Q_\tau$.  Then $a_\tau+c_\tau-2b_\tau$ is constant and $\detQ_\tau$ is a 
strictly decreasing function of $\tau$.  
Thus (\ref{eqn:secondrk2}) will hold for $Q$ if it holds for 
$Q_\tau$ for some $\tau>0$.  
In the same way we may increase both $b$ and $c$ so that $a+c-2b$ remains constant and
the determinant decreases, or we may increase $a$ and decrease $c$.  In this way we can find a path
$Q_\tau$ in the space of reduced matrices from any given $Q$ to a diagonal matrix
$\left(\begin{matrix} -1 & 0\\ 0 & -\tilde\detQ \end{matrix}\right),$ 
such that $a+b-2c$ is constant along the
path and the determinant decreases.  It follows that (\ref{eqn:secondrk2}) 
holds for $Q$, and  the 
inequality is strict unless 
$Q=\Diag.$

A similar but more involved argument establishes Conjecture \ref{conj:vec} for rank 3 forms.  
We briefly sketch the argument.
Let $Q$ be represented by a reduced matrix of rank 3 (see for example \cite{j}) and let $x_0\in\zz^3$.
By succesively adding $2\tau$ to a diagonal
entry and $\pm\tau$ to an off-diagonal entry one may find a path of reduced matrices 
from $Q$ to $\tilde Q$ along which
$\displaystyle\max_{x\equiv x_0 (2)}x^TQ x$ is constant and the absolute value of 
the determinant decreases.
One cannot always expect that $\tilde Q$ will be diagonal but one can show that the various matrices
which  arise all satisfy
$$\displaystyle\max_{x\equiv x_0 (2)}x^T\tilde Q x\ge -2-|\det \tilde Q|,$$
(with strict inequality unless $\tilde Q=\Diag$)
from which it follows that this inequality holds for all negative-definite rank 3 forms.

Finally note that for rank 2 forms, the determinant of the adjoint matrix $\ad Q$ is equal to
the determinant of $Q$.  Conjecture \ref{conj:covec} for rank 2 forms of odd
determinant now follows by applying
(\ref{eqn:firstrk2}) to $\ad Q$ and dividing by the determinant $\detQ$.
\endproof

\section{Determinant three}

In this section we describe to what extent we can generalise Elkies' proof of 
Theorem \ref{thm:elkies} to non-unimodular forms. 
For convenience we work with positive-definite forms.
We obtain the following result.

\begin{theorem}\label{thm:det3}
Let $Q$ be a positive-definite quadratic
form over the integers of rank $n$ and determinant $3$. Then 
either $Q$ has 
a characteristic vector $x$ with $Q(x,x)\le n+2$ or it has a characteristic covector $y$ with
$Q'(y,y)\le n-\frac{2}{3}$. Moreover, at least one of the above inequalities is strict unless 
$Q$ is diagonal. 
\end{theorem}

Given a  positive-definite integral quadratic form $Q$ of rank $n$, we consider lattices 
$L \subset L'$ in $\rr^n$ (equipped with the standard inner product), with $Q$ the intersection 
pairing of $L$, and $L'$ the dual lattice of $L$.  Note that the discriminant of the lattice $L$
is equal to the determinant of $Q$.

For any lattice $L \subset \rr^n$ and a vector $w \in \rr^n$ let $\theta_L^w$ be the generating 
function for the norms of vectors in $\frac{w}{2}+L$,
$$\theta_L^w(z)=\sum_{x\in L} e^{i\pi |x+\frac{w}{2}|^2 z};$$
this is a holomorphic function on the upper half-plane $H=\{z \,|\, \im(z)>0\}$.
The {\it theta series} of the lattice $L$ is $\theta_L=\theta_L^0$.

Recall that the modular group $\Gamma=\PSL_2(\zz)$ acts on $H$ and is generated by 
$S$ and $T$, where $S(z)=-\frac{1}{z}$ and $T(z)=z+1$. 
\begin{proposition}\label{prop:poisson}
Let $L$ be an integral lattice of odd discriminant $\detQ$, and $L'$ its dual lattice. Then 
\begin{eqnarray}
\theta_L(S(z))=\left(\frac{z}{i}\right)^{n/2}\detQ^{-1/2}\theta_{L'}(z) \label{p1} \\
\theta_L(TS(z))=\left(\frac{z}{i}\right)^{n/2}\detQ^{-1/2}\theta_{L'}^w(z) \label{p2} \\
\theta_{L'}(T^{\detQ}S(z))=\left(\frac{z}{i}\right)^{n/2}\detQ^{1/2}\theta_{L}^w(z) \label{p2'},
\end{eqnarray}
where $w$ is a characteristic vector in $L$.
\end{proposition}
\begin{remark}
Note that if $w \in L$ is a characteristic vector, then $\theta_{L'}^w$ is a generating function 
for the squares of characteristic covectors. Under the assumption that the discriminant of $L$ is odd, 
$\theta_L^w$ is a generating function for the squares of characteristic vectors.
\end{remark}
\proof
All the formulas follow from Poisson inversion \cite[Ch. VII, Proposition 15]{serre}.
We only need odd discriminant in (\ref{p2'}). Note that in $\theta_{L'}(z+\detQ)$ we can use
\begin{equation}\label{eqn:char}
\detQ |y|^2 \equiv |\detQ y|^2 \equiv (\detQ y,w) \equiv (y,w) \pmod 2
\end{equation}
and then apply Poisson inversion.
\endproof

\begin{corollary}\label{cor:sym}
Let $L_1$ and $L_2$ be integral lattices of the same rank and the same odd discriminant $\detQ$. 
Then 
$$R(z)=\frac{\theta_{L_1}(z)}{\theta_{L_2}(z)}$$
is invariant under $T^2$ and $ST^{2\detQ}S$. Moreover, $R^8$ is invariant under
$(T^2S)^{\detQ}$ and $ST^{\detQ-1}ST^{\detQ-1}S$.
\end{corollary}
\proof
Since $L$ is integral, $\theta_L(z+2)=\theta_L(z)$, hence $R$ is $T^2$ invariant. 
The squares of vectors in $L'$ belong to $\frac{1}{\detQ}\zz$, so
$\theta_{L'}(z+2\detQ)=\theta_{L'}(z)$. From (\ref{p1}) it follows that 
$R(S(z))=\frac{\theta_{L_1'}(z)}{\theta_{L_2'}(z)}$, which gives the $ST^{2\detQ}S$ invariance  
of $R$.

To derive the remaining symmetries of $R^8$ we need to use (\ref{p2}) and (\ref{p2'}).
Let $w$ be a characteristic vector in $L$. Clearly
$$\detQ |y+\frac{w}{2}|^2=\detQ |y|^2+ \detQ (y,w)+\frac{\detQ}{4}|w|^2$$
holds for any $y \in L'$, so it follows from (\ref{eqn:char}) that
$$\theta_{L'}^w(z+\detQ)=e^{i\pi \detQ |w|^2/4}\theta_{L'}^w(z).$$ 
Using (\ref{p2}) we now conclude
that $R^8$ is invariant under $TST^{\detQ}ST^{-1}=(ST^{-2})^{\detQ}$; the last equality follows 
from the relation $(ST)^3=1$ in the modular group. The remaining invariance of $R^8$
is derived in a similar way from (\ref{p2'}).
\endproof

From now on we restrict our attention to discriminant $\detQ=3$. Consider the subgroup $\Gamma_3$ 
of $\Gamma$ generated by $T^2$, $ST^6S$ and $ST^2ST^2S$. Clearly $\Gamma_3$ is a subgroup
of $\Gamma_+=\langle S,T^2 \rangle \subset \Gamma$. 

\begin{lemma}
A full set of coset representatives for $\Gamma_3$ in $\Gamma_+$ is $I,S,ST^2,ST^4$.
Hence a fundamental domain $D_3$ for the action of $\Gamma_3$ on the hyperbolic plane $H$ 
is the hyperbolic polygon with vertices $-1,-\frac{1}{3},-\frac{1}{5},0,1,i\infty$.
\end{lemma}
\proof
Call $x,y \in \Gamma_+$ {\it equivalent} if $y=zx$ for some $z \in \Gamma_3$. Let 
$x=T^{k_1}ST^{k_2}S \cdots T^{k_n}$ with all $k_i \ne 0$; then the {\it length} of $x$, $Sx$, 
$xS$ and $SxS$ is defined to be $n$. Any element $x \in \Gamma_+$ of length $n \ge 2$ 
is equivalent to one of the form $ST^{k}Sy$ with $k=0,\pm 2$ and length at most $n$. If 
$x=ST^{k}ST^{l}y$ with $k=\pm 2$ and length $n \ge 2$, then $x$ is equivalent to $ST^{l-k}y$,
which has length $\le n-1$. It follows by induction on length that any element of $\Gamma_+$
is equivalent to one with length at most $1$.
Moreover, if the element has length $1$, it is equivalent to $ST^k$, $k=2,4$.

Finally, recall that a fundamental domain for $\Gamma_+$ is $D_+=\{z\in H \,|\, -1\le \re(z)\le 1,\ 
|z|\ge 1\}$ so we can take $D_3$ to be the union of $D_+$ and $S(D_+ \cup T^2(D_+) \cup T^4(D_+))$. 
\endproof

\proof[Proof of Theorem \ref{thm:det3}]
Suppose that $L$ is a lattice of discriminant 3 and rank $n$ for which the square of any 
characteristic vector is at least $n+2$ and the square of any characteristic covector is
at least $n-\frac{2}{3}$. Let $\Diag$ be the lattice with intersection form
$(n-1)\langle1\rangle\oplus\langle3\rangle$;  recall from \cite{elkies} that $\theta_\Diag$
does not vanish on $H$. Then
$$R(z)= \frac{\theta_L(z)}{\theta_\Diag(z)}$$
is holomorphic on $H$ and it follows from Corollary \ref{cor:sym} that $R^8$ is invariant 
under $\Gamma_3$. We want to show that $R$ is bounded. We will use the following identities 
that follow from Proposition \ref{prop:poisson}:
$$R(S(z))=\frac{\theta_{L'}(z)}{\theta_{\Diag'}(z)},\ \ 
R(TS(z))=\frac{\theta_{L'}^w(z)}{\theta_{\Diag'}^w(z)},\ \
R(ST^{\detQ}S(z))= \frac{\theta_L^w(z)}{\theta_\Diag^w(z)}.
$$

Since the theta series of any lattice converges to $1$ as $z \to i\infty$, $R(z) \to 1$ as 
$z \to 0,i\infty$. By assumption the square of any characteristic covector for $L$ 
is at least as large as the square of the shortest characteristic covector for $\Diag$. Since the 
asymptotic behaviour as $z \to i\infty$ of the generating function for the squares of 
characteristic covectors is determined by the smallest square,
it follows from the middle expression for $R$ above that $R(z)$ is bounded as $z \to 1$. 
Similarly, using the condition on characteristic vectors and the right-most expression for
$R$ as $z \to i\infty$, it follows that $R(z)$ is bounded as $z \to -\frac{1}{3}$.
Note that $T^{-2}(1)=-1$ and $ST^6S(1)=-\frac{1}{5}$, 
so $R(z)$ is also bounded as $z \to -1,-\frac{1}{5}$. 

Let $f$ be the function on $\Sigma=H/\Gamma_3$ induced by $R^8$. Then $f$ is holomorphic and 
bounded, so it extends to a holomorphic function on the compactification of $\Sigma$. It 
follows that $R(z)=1$, so the theta series of $L$ is equal 
to the theta series of $\Diag$. Then $L$ contains $n-1$ pairwise orthogonal vectors of square
$1$, so its intersection form is $(n-1)\langle1\rangle\oplus\langle3\rangle$.
\endproof

\section{Applications}
\label{sec:examples}
In this section we consider applications to rational homology spheres
and four-ball genus of knots.  
We begin with the proof of Theorem \ref{cor:det3}.

\proof[Proof of Theorem \ref{cor:det3}.]
Suppose that $Y=\partial X$ and that $Q$ is the intersection form on $H_2(X;\zz)$.  
Then $Q$ is a quadratic form of determinant $\pm3$.
For any $\spincs\in\spinc(X)$, let $c(\spincs)$ denote the image of the first Chern class 
$c_1(\spincs)$ modulo torsion. Then $c(\spincs)$  is a characteristic
covector for $Q$; moreover if $\spincs|_Y$ is \spin then $c(\spincs)$ is $Qx$ for some
characteristic vector $x$.
The result now follows from Theorem \ref{thm:det3} and (\ref{eqn:thm9.6}).\endproof

Conjectures \ref{conj:vec} and \ref{conj:covec} imply the following more general statement.

\begin{conjecture}
\label{conj:app}
Let $Y$ be a rational homology sphere with $|H_1(Y;\zz)|=h$.
If Y bounds a negative-definite four-manifold $X$ with no torsion in $H_1(X;\zz)$ then
$$\displaystyle\min_{\spinct_0\in\spin(Y)}d(Y,\spinct_0)\ge (1-h)/4,$$
and
$$\displaystyle\max_{\spinct\in\spinc(Y)}d(Y,\spinct)\ge
\left\{\begin{array}{ll}
\left(\displaystyle1-\frac1h\right)/4 & \mbox{if $h$ is odd,}\\[5mm]
1/4 & \mbox{if $h$ is even.}
\end{array}\right.$$
If equality holds in either inequality the intersection form of $X$ is $\Diag_h$.
\end{conjecture}

More generally if $Y$ bounds $X$ with torsion in $H_1(X;\zz)$,
the absolute value of the determinant of the intersection pairing of $X$ divides $h$ 
with quotient a square (see for example \cite[Lemma 2.1]{bs}).
One may then deduce inequalities as above corresponding to each
choice of determinant; care must be taken since for example not
all \spin structures on $Y$ extend to \spinc structures on $X$.

\begin{remark}
Given a rational homology sphere $Y$ bounding $X$ with no torsion in $H_1(X;\zz)$, 
the intersection pairing of $X$ gives a presentation matrix for $H^2(Y;\zz)$ (and also
determines the linking pairing of $Y$).
There should be analogues of Conjectures \ref{conj:vec} and \ref{conj:covec}
which restrict to forms presenting a given group (and inducing a given linking pairing).  
These should give stronger bounds than those in Conjecture \ref{conj:app}.
\end{remark}

\subsection{Seifert fibred examples}
\label{subsec:sfs}

In Examples \ref{ex:cantbound1} and \ref{ex:cantbound2}  we list families 
of Seifert fibred spaces $Y$ which bound positive-definite but not 
negative-definite four-manifolds. It follows as in \cite[Theorem 10.2]{fs} that
no multiple of $Y$ bounds a negative-definite four-manifold. 
In Examples \ref{ex:diag1}
through \ref{ex:diag3} we list families of Seifert fibred spaces which can only bound the 
diagonal negative-definite form $\Diag_\detQ$ (or sometimes $\Diag_1$).  
We found these examples using predictions
based on Conjecture \ref{conj:app} and verified them using Donaldson's theorem via
Proposition \ref{prop:cantbound}. Finally, in Example \ref{ex:diag4} we exhibit a family of Seifert
fibred spaces which according to the conjecture  can only bound $\Diag_\detQ$. For this
family the method of Proposition \ref{prop:cantbound} does not apply.

In what follows we extend the definition of $\Diag_1$ to include the trivial form on the 
trivial lattice.  Also note that a lattice uniquely determines a quadratic form, and
a form determines an equivalence class of lattices; in the rest of this section we use the terms
lattice and form interchangeably.

\begin{definition}
Let $L$ be a  lattice  of rank $m$ and determinant $\detQ$.
We say $L$ is {\em \rigid} if any embedding of $L$ in $\zz^n$
is contained in a $\zz^m$ sublattice.  We say $L$ is {\em almost-\rigid}
if any embedding of $L$ in $\zz^n$ is either contained in a
$\zz^m$ sublattice, or contained in a $\zz^{m+1}$ sublattice
with orthogonal complement spanned by a vector $v$ with $|v|^2=\detQ$.
\end{definition}

\begin{proposition}
\label{prop:cantbound}
Let $Y$ be a rational homology sphere and let $h$ be the order of $H_1(Y;\zz)$. 
Suppose $Y$ bounds a positive-definite four-manifold $X_1$ with $H_1(X_1;\zz)=0$.
Let $Q_1$ be the intersection pairing of $X_1$ and let $m$ denote its rank.

If $Q_1$ does not embed into $\zz^n$ for any $n$ 
then $Y$ cannot bound a negative-definite four-manifold.

If $Q_1$ is rigid and $Y$ bounds a negative-definite $X_2$ then $h$ is a square and 
$Q_2=\Diag_1$; if $h>1$, then there is torsion in $H_1(X_2;\zz)$.

If $Q_1$ is almost-rigid and $Y$ bounds a negative-definite $X_2$ then either
\begin{itemize}
\item $Q_2=\Diag_h$ or
\item $Q_1$ embeds in $\zz^m$,  $h$ is a square and 
$Q_2=\Diag_1$; if $h>1$, then there is torsion in $H_1(X_2;\zz)$.
\end{itemize}
\end{proposition}

\proof
Suppose $Y$ bounds a negative-definite $X_2$ with intersection pairing $Q_2$.
Then $X=X_1\cup_Y -X_2$ is a closed positive-definite
manifold.  The Mayer-Vietoris sequence for homology and Donaldson's theorem
yield an embedding $\iota:Q_1\oplus -Q_2\to\zz^{m+k}$, where $k$ is the rank of $Q_2$.

If the image of $Q_1$ under $\iota$ is contained in a $\zz^m$ sublattice, then the image of $-Q_2$ is
contained in the orthogonal $\zz^k$ sublattice.  Now consider the Mayer-Vietoris sequence for cohomology:
$$\begin{array}{cccccccc}
0\longrightarrow& H^2(X;\zz) & \longrightarrow & H^2(X_1;\zz)&\!\!\!\oplus\!\!\!&H^2(X_2;\zz) 
&\longrightarrow&H^2(Y;\zz),\\
&&&\|&&\|&&\\
&&&Q_1'&&\!\!\!\!\!-Q_2'\oplus T_2&&
\end{array}$$
where $T_2$ is the torsion subgroup and $Q'$ denotes the dual lattice to $Q$.
This yields an embedding $\iota':\zz^{m+k}\to Q_1'\oplus -Q_2'$.  The mapping $\iota'$ 
is hom-dual to $\iota$ and hence also decomposes orthogonally, sending
$\zz^m$ to $Q_1'$ and $\zz^k$ to $-Q_2'$.  The image of $\zz^m$ in  $Q_1'$ has index $\sqrt{h}$,
since $h$ is the determinant of $Q_1$.  (In general if $L_1\subset L_2$ are lattices of the same
rank then the square of the index $[L_2:L_1]$ is the quotient of their discriminants.)
The restriction map from $H^2(X_1;\zz)$ to $H^2(Y;\zz)$ is onto, so its kernel $K$ is a subgroup
of $\zz^m$ of index $\sqrt h$. It follows that $\zz^m/K$ injects into $T_2$ and that the image
of $T_2$ in $H^2(Y;\zz)$ has order $t \ge \sqrt h$. Then by \cite[Lemma 2.1]{bs}, 
$t=\sqrt h$ and $Q_2$ is unimodular.
Since $-Q_2$ is a sublattice of $\zz^k$ we have $Q_2=\Diag_1$.

Suppose now that the image of $Q_1$ under $\iota$ is contained in a $\zz^{m+1}$ sublattice, and its
orthogonal complement in $\zz^{m+1}$ is spanned by a vector $v$ with $|v|^2=h$.  
Then the image of $-Q_2$ is a sublattice of $(k-1)\langle1\rangle\oplus \langle h\rangle$; it therefore has determinant
at least $h$.  On the other hand its determinant divides $h$ \cite[Lemma 2.1]{bs}.  It follows
that $Q_2$ is equal to $\Diag_h$.
\endproof

If $Y$ is the Seifert fibred space
$Y(e;(\alpha_1,\beta_1),(\alpha_2,\beta_2),(\alpha_3,\beta_3))$, let 
$$k(Y)=e\alpha_1\alpha_2\alpha_3+\beta_1\alpha_2\alpha_3+
\alpha_1\beta_2\alpha_3+\alpha_1\alpha_2\beta_3.$$
If $k(Y)\ne 0$ then $Y$ is a rational homology sphere and
$|k(Y)|$ is the order of $H_1(Y;\zz)$. Furthermore, if $k(Y)<0$ then $Y$ bounds a
positive-definite plumbing. 
For our conventions for lens spaces and Seifert fibred spaces see \cite{bs}. Recall in particular
that $(\alpha_i,\beta_i)$ are coprime pairs of integers with $\alpha_i\ge2$.  We will also assume here
that $1\le\beta_i<\alpha_i$.

\begin{example}
\label{ex:cantbound1}
Seifert fibred spaces $Y=Y(-2;(\alpha_1,\beta_1),(\alpha_2,\beta_2),(\alpha_3,\beta_3))$ with
$$\frac{\alpha_1}{\beta_1}\le2,\quad\frac{\alpha_2}{\beta_2},\frac{\alpha_3}{\beta_3}<2,\quad k(Y)<0,$$
cannot bound negative-definite four-manifolds.
\end{example}

\begin{figure}[htbp]
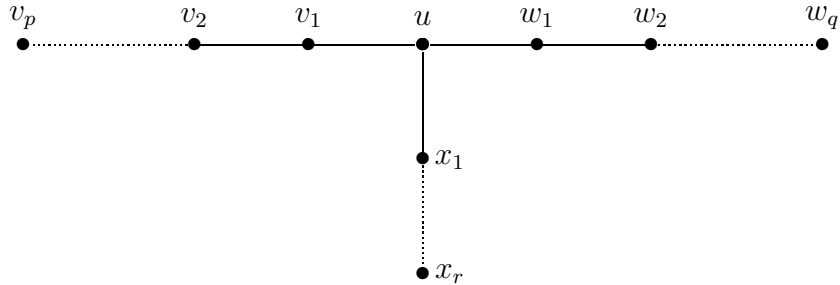

  \begin{center}
\ifpic
    \leavevmode
    \[ \xy/r1.8pc/:
  (0,0)="centre", "centre"+/u2ex/*{u},
  (-2,0)="3,1", "3,1"+/u2ex/*{v_1},
  (-4,0)="3,2", "3,2"+/u2ex/*{v_2},
  (-7,0)="3,s", "3,s"+/u2ex/*{v_p},
  (2,0)="1,1", "1,1"+/u2ex/*{w_1},
  (4,0)="1,2", "1,2"+/u2ex/*{w_2},
  (7,0)="1,s", "1,s"+/u2ex/*{w_q},
  (0,-2)="2,1", "2,1"+/r2ex/*{x_1},
  (0,-4)="2,s", "2,s"+/r2ex/*{x_r},
  "centre"*{\bullet};      
  "1,2" **\dir{-},         
  "1,1"*{\bullet};         
  "1,2"*{\bullet};         
  "1,s" **\dir{.},         
  "1,s"*{\bullet};         
  "centre"*{\bullet};      
  "2,1" **\dir{-},         
  "2,1"*{\bullet};         
  "2,s" **\dir{.},         
  "2,s"*{\bullet};         
  "centre"*{\bullet};      
  "3,2" **\dir{-},         
  "3,1"*{\bullet};         
  "3,2"*{\bullet};         
  "3,s" **\dir{.},         
  "3,s"*{\bullet};         
\endxy \]
\else \vskip 5cm \fi
   \begin{narrow}{0.3in}{0.3in}
    \caption{{\bf Plumbing graph.}
        }
    \label{fig:plumbing}
   \end{narrow}
  \end{center}
\end{figure}

\proof
Note that $Y$ is the boundary of the positive-definite plumbing shown in Figure \ref{fig:plumbing}, where
vertices $u$, $v_1$, $w_1$ and $x_1$ have square 2 and $v_2$ and $w_2$ have square at least 2.
This lattice does not admit an embedding in any $\zz^n$. To see this let $e_1,\ldots,e_n$ be 
the standard basis of $\zz^n$. 
The vertex $u$ must map to  an element of square $2$, which we may suppose
is $e_1+e_2$.  The 3 adjacent vertices must be mapped to elements of the form
$e_1+e_3$, $e_1-e_3$ and $e_2+e_4$.  Now we see that it is not possible to map the remaining
2 vertices $v_2$ and $w_2$; we are only able to further extend the map along the leg of
the graph emanating from the vertex mapped to $e_2+e_4$.
\endproof

\begin{example}
\label{ex:cantbound2}
Seifert fibred spaces $Y=Y(-2;(\alpha_1,\beta_1),(\alpha_2,\alpha_2-1),(\alpha_3,\alpha_3-1))$ with
$$\alpha_2,\alpha_3\ge\frac{\alpha_1}{\beta_1},\quad \alpha_3\ge 3,\quad k(Y)<0,$$
cannot bound negative-definite four-manifolds unless 
$$\beta_1=1,\quad \min(\alpha_2,\alpha_3)=\alpha_1.$$
In the latter case, if $Y$ bounds a negative-definite $X$ then the intersection pairing of $X$ is
$\Diag_1$ and the torsion subgroup of $H_1(X;\zz)$ is nontrivial.
\end{example}
\proof
In this case $Y$ is again the boundary of a positive-definite plumbing as in Figure \ref{fig:plumbing}.
The vertices $u$, $v_i$ and $w_j$ have square 2, and $p=\alpha_2-1$, $q=\alpha_3-1$.
Vertex $x_1$ has square $a=\lceil\frac{\alpha_1}{\beta_1}\rceil$.  If 
$\frac{\alpha_1}{\beta_1}=\min(\alpha_2,\alpha_3)=a$ then by inspection this pairing is \rigid
with determinant $a^2>1$;
otherwise it does not admit any embedding into $\zz^n$.  For more details see the proof of
Example \ref{ex:diag2}.
\endproof

\begin{example}
\label{ex:diag1}
The only negative-definite pairing that $L(p,1)$ can bound is the diagonal form $\Diag_p$ unless
$p=4$ in which case it may also bound $\Diag_1$.  (Note that
$L(p,1)$ is the boundary of the disk bundle over $S^2$ with intersection pairing $\langle-p\rangle$.)
\end{example}
\proof
Observe that $L(p,1)$ is the boundary of the positive-definite plumbing $A_{p-1}$. If $p\ne 4$ then
up to automorphisms of $\zz^n$ there is a unique embedding of $A_{p-1}$ in $\zz^n$; the image is
contained in a $\zz^p$ and its orthogonal complement in $\zz^p$ is generated by the vector 
$(1,1,\ldots,1)$.  Hence $A_{p-1}$ is almost-\rigid and does not embed in $\zz^{p-1}$.
However, $A_3$ also admits an embedding in $\zz^3$.
\endproof

\begin{example}
\label{ex:diag2}
If $Y=Y(-2;(\alpha_2\beta_1+1,\beta_1),(\alpha_2,\alpha_2-1),(\alpha_3,\alpha_3-1))$ with
$\alpha_3>\alpha_2,$
then the only negative-definite pairing that $Y$ may bound is the diagonal form $\Diag_{|k(Y)|}$
unless 
$$\beta_1=1,\quad \alpha_3=\alpha_2+1.$$
In the latter case the only negative-definite pairings that $Y$ may bound are  $\Diag_{|k(Y)|}$
and $\Diag_1$.
\end{example}
\proof
Note this is a borderline case of Example \ref{ex:cantbound2}.  In the notation of that example 
$\alpha_2=a-1$. The positive-definite plumbing is similar to that in Example \ref{ex:cantbound2}
with $r=\beta_1$; also the vertices $x_l$ with $l>1$ all have square 2.  Denote the
pairing associated to this plumbing by $Q$.
We consider an embedding of $Q$ into $\zz^n$. Let $e_i$, $f_j$ and $g_l$ denote 
unit vectors in $\zz^n$.    Without loss of generality 
the vertex $u$ maps to $e_1+f_1$.  Then $v_i$ maps to 
$e_{i-1}+e_i$ and $w_j$ maps to $f_{j-1}+f_j$.  

Now consider the image of $x_1$.  This may map
to $e_1-e_2+\cdots\pm e_{a-1}+g_1$; then $x_l$ maps to $g_{l-1}+g_l$ for $l>1$.  Thus the image of 
$Q$ is contained in a $\zz^{p+q+r+2}$ sublattice.  The determinant of $Q$ is 
$|k(Y)|=\alpha_2{}\!^2\beta_1+\alpha_2+\alpha_3$ (note $k(Y)<0$).  
The orthogonal complement of $Q$ in $\zz^{p+q+r+2}$
is spanned by the vector $\sum(-1)^{i-1}e_i+\sum(-1)^j f_j+\alpha_2\sum(-1)^l g_l$, whose square is
$|k(Y)|$.  Up to automorphism this is the only embedding of $Q$ into $\zz^n$ unless
$\alpha_3=a$ and $\beta_1=1$.  In this case $x_1$ may map to the alternating sum $f_1-f_2+\cdots\pm f_a$;
the image of the resulting embedding is contained in $\zz^{p+q+r+1}$.
\endproof

\begin{example}
\label{ex:diag3}
If $Y=Y(-1;(3,1),(3a+1,a),(5b+3,2b+1))$ with $k(Y)<0$,
then  the only negative-definite pairing that $Y$ 
may bound is the diagonal form $\Diag_{|k(Y)|}$ unless $a=b=1$ in which 
case it may also bound $\Diag_1$. 
\end{example}
\proof
Note that the condition $k(Y)<0$ implies 
$a=1$ or $b=0$ or $a=b+1=2$.
Again, $Y$ is the boundary of a positive-definite plumbing as in Figure \ref{fig:plumbing},
with $p=a$, $q=b+1$ and $r=1$.
The vertex $u$ has square 1, $w_1$ and $x_1$ have square 3, $v_1$ has square 4.
If $a>1$ then $v_j$ has square 2 for $j>1$. If $b>0$ then $w_2$ has square 3, and any remaining 
$w_i$ has square 2. 
Denote the pairing associated to this plumbing by $Q$.
We consider an embedding of $Q$ into $\zz^n$. Let $e_i$ denote 
unit vectors in $\zz^n$.    Without loss of generality 
the vertex $u$ maps to $e_1$, $x_1$ maps to $e_1+e_2+e_3$ and $w_1$ maps to $e_1-e_2+e_4$.
Then $v_1$ has to map to $e_1-e_3-e_4+e_5$. Now $w_2$, if present, has to map to 
$e_4+e_5+e_6$ or $-e_2+e_3+e_5$; the second possibility only works if $a=b=1$. 
Finally $v_2$, if present, has to map to $e_5-e_6$. 
The reader may verify that $Q$ is almost-\rigidpunc. 
\endproof

\begin{example}
\label{ex:diag4}
Let $Y_a=Y(-2;(2,1),(3,2),(a,a-1))$ with $a\ge 7$. Then $h=k(Y)=a-6$,
$$\displaystyle\min_{\spinct_0\in\spin(Y)}d(Y,\spinct_0)= (1-h)/4$$
and
$$\displaystyle\max_{\spinct\in\spinc(Y)}d(Y,\spinct)=
\left\{\begin{array}{ll}
\left(\displaystyle1-\frac1h\right)/4 & \mbox{if $h$ is odd,}\\[5mm]
1/4 & \mbox{if $h$ is even.}
\end{array}\right.
$$
If $a$ is $7$ or $9$ then the only negative-definite form $Y_a$ bounds is $\Diag_h$. If
Conjecture \ref{conj:app} holds then the same is true for all $Y_a$. 
\end{example}
\proof
$Y_a$ is the boundary of the negative-definite plumbing with intersection pairing given by
$$Q=\left(\begin{matrix} -1 & 1 & 1 & 1 \\ 1 & -2 & 0 & 0 \\ 
1 & 0 & -3 & 0 \\ 1 & 0 & 0 & -a \end{matrix}\right), $$
which represents $3\langle-1\rangle \oplus\langle-a+6\rangle$. The computations of $d(Y)$ 
follow as in \cite{os6}. The claim for $Y_7$ follows from the discussion following Theorem
\ref{thm:elkies}. The claim for $Y_9$ follows from Theorem \ref{cor:det3}.
\endproof

\subsection{Four-ball genus of Montesinos knots}
\label{subsec:knots}
Let $K$ be a knot in $S^3$ and let $g$ denote its Seifert genus. 
The four-ball genus $g^*$ of $K$ is the minimal genus of a 
smooth surface in $B^4$ with boundary $K$.  A classical result of Murasugi states that
$g^*\ge|\sigma|/2$, where $\sigma$ is the signature of $K$.  If this lower bound is
attained then the double branched cover of $S^3$ along $K$ bounds a definite four-manifold
with signature $\sigma$.
The double branched cover of the Montesinos knot or link 
$M(e;(\alpha_1,\beta_1),(\alpha_2,\beta_2),(\alpha_3,\beta_3))$ is
$Y(-e;(\alpha_1,\beta_1),(\alpha_2,\beta_2),(\alpha_3,\beta_3))$.
(For more details see \cite{bs}.)

The following generalises an example of Fintushel and Stern \cite{fs}.
\begin{example}
\label{ex:knots1}
The pretzel knot $K(p,-q,-r)=M(2;(p,1),(q,q-1),(r,r-1))$ for odd $p$, $q$ and $r$ satisfying
$$q,r>p>0\quad \mbox{and}\quad pq+pr-qr\ \mbox{is a square}$$
is algebraically slice but 
has $g^*=1$.
\end{example}
\proof
The knot has a genus 1 Seifert surface yielding the Seifert matrix
$$M=\left(\begin{matrix}
\frac{p-r}2 & \frac{p+1}2\\
\frac{p-1}2 & \frac{p-q}2
\end{matrix}\right).$$
The vector $x=(p-l,r-p)$, where $l=\sqrt{pq+pr-qr}$, satisfies $x^TMx=0$,
demonstrating the knot is algebraically slice.
The double branched cover $Y$ of the knot has $k(Y)=-l^2$.
From Example \ref{ex:cantbound2} we see that $Y$  does not 
bound a rational homology ball.  It follows that $0<g^*\le g = 1$. 

It is shown by Livingston \cite{l} that  $K(p,-q,-r)$ has $\tau=1$, where $\tau$ is the
\ozsvath-\szabo\ knot concordance invariant. This also gives $g^*=1$.
\endproof

\begin{example}
\label{ex:knots2}
The Montesinos knot $K_{q,r}=M(2;(qr-1,q),(r+1,r),(r+1,r))$ with
odd $q\ge 3$ and even $r\ge 2,$ has signature $\sigma=1-q$ and has
$$g=g^*= \frac{q+1}2.$$
Computations suggest that the Taylor invariant of $K_{q,r}$ is $\frac{q-1}2$.
\end{example}
\proof
The knot $K_{q,r}$ is equal to $M(0;(qr-1,q),(r+1,-1),(r+1,-1))$. It is easily seen
that $K_{q,r}$ has a spanning surface with genus $\frac{q+1}2$. Using the resulting Seifert
matrix one gets the formula for the signature. The 
double branched cover $Y$ of $K_{q,r}$ has $k(Y)<0$.
From Example \ref{ex:cantbound2} we see that $Y$  does not
bound a negative-definite four-manifold; the genus formula follows.

We have computed the Taylor invariant of $K_{q,r}$ for $q< 10000$ and any $r$. 
\endproof

\begin{remark}
We have discussed Conjectures \ref{conj:vec} and \ref{conj:covec} with Noam Elkies.  He has
suggested an alternative proof of Theorem \ref{thm:det3} using gluing of lattices 
\cite{elkies2}.
His proof works for odd determinants $\detQ$ up to 11, under the additional assumption 
that there is an element of $L'$ whose square is congruent to $1/\detQ$ modulo 1.
He also indicated a way to remove this assumption. 

Using his result it follows that $\Diag_{a-6}$ is the only negative-definite form bounded by the manifold
$Y_a$ of Example \ref{ex:diag4} for $a=11$,
$13$, $15$ and $17$.
\end{remark}

\end{document}